\documentclass[11pt,a4paper]{amsart}
\usepackage[left=3.8cm, right=3.3cm, top=4.0cm, bottom=4.0cm]{geometry}
\usepackage{color}
\usepackage[latin1]{inputenc}
\usepackage{amssymb}
\usepackage{amsthm}
\theoremstyle{plain}
\usepackage{epsfig}
\usepackage{url}
\usepackage{hyperref}
\usepackage{graphicx}

\newtheorem{thm}{Theorem}[section]
\newtheorem*{main}{Main result}
\newtheorem{cor}[thm]{Corollary}
\newtheorem{lem}[thm]{Lemma}
\newtheorem{pro}[thm]{Proposition}

\theoremstyle{definition}
\newtheorem*{rem}{Remark}
\newtheorem*{rems}{Remarks}
\newtheorem{exa}[thm]{Example}
\newtheorem{exas}[thm]{Examples}

\newcommand\guedel[1]{\marginpar{GR del}}
\newcommand\ganzweg[1]{}
\newcommand{\beq}{\begin{equation}}
\newcommand{\eeq}{\end{equation}}
\newcommand{\R}{\mathbb{R}}

\begin{document}
\title[Ricci almost solitons: complete examples]
{Ricci almost solitons: complete examples}
\author[W.K\"uhnel \& H.B.Rademacher]{Wolfgang K\"uhnel and Hans-Bert Rademacher}
\date{2026-02-11, revised version}

\thanks{}
\begin{abstract}
  {\scriptsize We construct new examples of various solitons as warped products.
There are classes of complete Ricci almost solitons and complete    
Ricci-Bourguignon solitons that can be explicitly described in terms of
elementary functions.}
\end{abstract}
\maketitle

\section{Introduction and results}
A {\it Ricci soliton} is a self-similar solution of the Ricci flow, i.e.,
a Riemannian metric $g$ satisfying the equation
$${\rm Ric} + \frac{1}{2}\mathcal{L}_Xg = \lambda g $$
for the Lie derivative $\mathcal{L}_Xg$ in the direction of
some vector field $X$ and a constant $\lambda$.
 Any Einstein space is
 trivially a Ricci soliton with $X=0$.
 If in addition $X$ is assumed to be a non-trivial conformal vector field,
 then it immediately follows that the Ricci soliton is Einstein.
 This case is often discussed in the literature, see \cite{JW}, \cite{D}.
 Besides the standard sphere, compact Ricci solitons seem to be rare
 \cite{De1, De2}.
 Sophisticated analytic and topological
 methods have been employed for obtaining non-compact Ricci solitons
 \cite{Wi}, \cite{CRZ}, for the approach by doubly warped products see
 \cite{Iv}.

 \medskip
  If $X$ is the gradient of a function $f$ then it is called
  a {\it gradient Ricci soliton}, characterized by the equation
  \begin{equation}
  	\label{eq:ricci}  
  {\rm Ric} + \nabla^2 f = \lambda g \,,
\end{equation}
  where $\nabla^2f$ is the Hessian.
  If $\lambda$ is a non-constant function then one talks about an
  {\it Ricci almost soliton} or {\it gradient Ricci almost soliton},
  respectively.
  Einstein spaces that are Ricci almost solitons with a nontrivial $X$
  have to admit a nontrivial conformal vector field. These spaces
  are classified in a finite list \cite{K}, also in
  pseudo-Riemannian geometry \cite{KR1}.

  \medskip
  A family $t \mapsto g(t)$ of Riemannian metrics
satisfying the evolution equation
  \begin{equation*}
  	\frac{\partial g}{\partial t}=
  	-2 \left(\textrm{Ric}-\rho S g\right)  	
  \end{equation*}
for a constant $\rho \in \mathbb{R}$
is said to 
  evolve by the \emph{Ricci-Bourguignon flow}.
  Here $S$ is the scalar curvature of $g=g(t).$
  We will consider two particular cases.
  For $\rho=0$ this is the Ricci flow,
  for $\rho=1/n$ the right hand side equals
  the traceless Ricci tensor.
  This flow was introduced by Bourguignon~\cite[3.24]{B}.
 
  \medskip
  
  Then a 
  \emph{gradient Ricci-Bourguignon soliton} 
  is a Riemannian metric $g$ together with
  a smooth function $f$ such that
  the following equation holds for some
  constant $\lambda:$
  \begin{equation}
  	\label{eq:ricci-bourgouignon}
  \textrm{Ric} + \nabla^2 f-
  \rho S g= \lambda g.
\end{equation}
Gradient Ricci-Bourguignon 
solitons are also named
\emph{gradient $\rho$-Einstein solitons,}
if $\lambda$ in Equation~\eqref{eq:ricci-bourgouignon}
is allowed to be a function,
then the metric is also called
\emph{gradient Ricci-Bourguignon almost soliton}.
For $\lambda=0$ the gradient Ricci-Bourguignon (almost) soliton
is called \emph{steady.}

  \medskip 
  In the present article we focus on complete and non-compact solutions
  of Equation~\eqref{eq:ricci} and
  Equation~\eqref{eq:ricci-bourgouignon}
  on warped products. We find new Ricci almost solitons of this type in Section~\ref{sec:complete},
  and in Theorem~\ref{thm:Ricci-Bourguignon} we classify for
  $n \ge 3$ and $\rho=1/n$ all complete
  and non-compact steady Ricci-Bourguignon solitons
  of this type, which then have a Ricci flat fibre. All solutions are explicit
  in terms of elementary functions. Ricci-Bourguignon solitons
  as generalized warped products were investigated in \cite{KO}.
 \begin{main}
All complete warped products with a $1$-dimensional basis
satifying the equation
$\textrm{Ric} + \nabla^2 f-  \frac{S}{n} g = 0$
of a steady gradient Ricci-Bourguignon soliton can
be classified by a family of explicitly given functions $h$ and $f$
depending on finitely many parameters.
\end{main}
      
  Throughout this paper ${\rm Ric}$ and
  $({\rm Ric})^\circ =  {\rm Ric} - \frac{S}{n} g$
  will denote the Ricci tensor and the traceless Ricci tensor, respectively,
  $S$ is the scalar curvature.
  
  \section{Basic formulas}
We consider an $n$-dimensional Riemannian manifold
$(M,g) \cong \mathbb{R} \times (M_*,g)$ with $n \geq 2$ and
with a warped product metric
$$ds^2 = dt^2 + h'^2(t)ds_*^2$$
and a function $f(t)$ defined on the basis.
Unless stated otherwise, it will be assumed that $h$ and $f$ are
smooth functions and $h'$ is not constant and non-zero.
Isolated zeros of $h'$ are allowed if the metric singularity is
removable, as it is the case for polar coordinates.
The gradients of $h$ and $f$ appear as
$h'\partial_t$ and $f'\partial_t$, respectively, and the
Hessian tensors are
$\nabla^2 h=h'' g$ on the one hand and
\begin{equation*}
\nabla_{\partial_t}\nabla f=
f'' \partial_t
\end{equation*}
and
\begin{equation*}
\nabla_X\nabla f=\nabla_X\left(f'\partial_t\right)=
f'\nabla_{X}\partial_t=
\frac{f'}{h'}\nabla_X\nabla h=
\frac{f' h''}{h'}X
\end{equation*} 
for any $X$ perpencidular to $\partial_t$ on the other hand.

\medskip
The Ricci tensors of $g$ and $g_*$ will be denoted by Ric and ${\rm Ric}_*,$
respectively.
Any equation of the type 
\begin{equation}\label{Ausgangsgleichung}
a \textrm{Ric}+b \nabla^2 f+ c g=0
\end{equation}
with functions $a, b, c$ on $M$ leads to equations of the form
\begin{eqnarray*}
a \textrm{Ric}(\partial_t,\partial_t)+
b f''+c &=&0\\
a \textrm{Ric}(X,X)+
b \frac{f' h''}{h'} g(X,X)+c g(X,X)&=&0\,.
\end{eqnarray*}
For the warped product metric above we obtain in direction $\partial_t:$
\begin{equation}
\label{eq:t-direction}
-(n-1)a\frac{h'''}{h'}+
b f''+c=0
\end{equation}
and in any direction $X$ perpendicular to 
$\partial_t:$
\begin{equation}
\label{eq:orthogonal}
a \textrm{Ric}_*
+
\left\{-a (n-2)h''^2 -
a h' h'''
+b f'h'h'' + c h'^2
\right\}g_*
=0.
\end{equation}
From Equation~\eqref{eq:orthogonal} we see that $(M_*,g_*)$ is Einstein,
so for $n-1 \geq 3$ there is a constant 
$\lambda_*$ (the normalized Einstein constant of $M_*$) satisfying
\begin{equation}\label{Hauptgleichung}
a (n-2)h''^2+
a h' h'''
-b f'h'h''-c h'^2
=\lambda_* a.
\end{equation}
For given constants $a, b$ and given $h'$ one can easily determine
an appropriate $f'$ and a function $c$ such that the equations
above are satisfied.

\medskip
We eliminate $c$ from the equations by inserting
Equation~\eqref{eq:t-direction} into Equation~\eqref{Hauptgleichung}:
\begin{equation}
a (n-2)h''^2 +
a h' h'''
- b f'h'h''
- (n-1)ah'h''' + bf''h'^2 
= \lambda_* a
\end{equation}
or, equivalently,
\begin{equation}\label{alteGleichung}
a (n-2)\big[h''^2 - h'h'''\big] 
 + bh'\big[f''h' - f'h'' \big] 
= \lambda_* a
\end{equation}
or
\begin{equation}\label{neueGleichung}
a (n-2)h''^2 \Big( \frac{h'}{h''}\Big)' 
 + bh'^3 \Big( \frac{f'}{h'}\Big)' 
= \lambda_* a.
\end{equation}

Therefore one can express $f'$ as $h'$ times an integral of
$\frac{a}{b} (n-2)\frac{h''^2 - h'h'''-\lambda_*}{h'^3}$
and determine $f$ itself by an additional integration step,
compare \cite[Ex.1.5]{PRRS}.
Finally, $c$ is obtained by the trace of $a{\rm Ric} + b\nabla^2 f$.
However, in general $c$ won't be constant, so this method does not
produce Ricci solitons but only Ricci almost solitons.
Considering only the traceless part leads to the following:

\begin{pro} \label{given-h}
  For any given Einstein space $M_*$ and any smooth function $h'$
  there is a function $f$ such that the warped product
  $\mathbb{R} \times_{h'} M_*$ admits at least locally a solution of
  $(a{\rm  Ric} + b\nabla^2 f)^\circ = 0$.
  For $a=b=1$ the left hand side is nothing but the traceless
  Bakry-\'Emery-Ricci tensor
  \cite{WW}.

  If $h'$ is globally defined and everywhere strictly positive, then
  the warped product   $\mathbb{R} \times_{h'} M_*$ is complete whenever
  $M_*$ is complete. This leads to complete solutions of
    $({\rm  Ric} + \nabla^2 f)^\circ = 0$ (Ricci almost solitons).
  \end{pro}

{\it Proof.} Equation~\eqref{alteGleichung} expresses the fact that
both Equation~\eqref{eq:t-direction} and Equation~\eqref{eq:orthogonal}
are simultaneously satisfied with a certain function $c(t)$.
This implies Equation~\eqref{Ausgangsgleichung} for this unknown $c$
or, equivalently, $(a {\rm  Ric} + b\nabla^2 f)^\circ = 0$.
The second claim is obvious from the explicit formula for $f'$ in
terms of $h' > 0$. \hfill $\Box$

\begin{pro}  {\rm (special cases)}\label{specialcases}
  \begin{enumerate}
\item The case $a=1$ and $b=0$ and, consequently, 
$(n-2)\big[h''^2 - h'h'''\big] = \lambda_*$ leads to an Einstein space
and the standard (non-constant) solutions
$h'(t) = \sin t, \cos t, \sinh t, \cosh t, e^t, t$ of $(\nabla^2h)^\circ = 0$
on such \cite{K}, \cite{KR1}.
\item
Essentially the same holds under the
assumption that $a=b=1$ and $M$ is Einstein, compare \cite[Thm.1.3]{PRRS}.
Up to scaling, the list
  of Einstein spaces admitting a function $f$ satisfying
  $\nabla^2 f = \lambda g$
  for some function $\lambda$ consists of warped products
  $dt^2 + (f'(t))^2ds_*^2$ with a $1$-dimensional basis where $f(t)$
  is a solution of the oscillator equation $f'' + \sigma f = 0$ where
  $\sigma$ is a normalized scalar curvature ($\sigma = 1$ for the unit sphere)
  \cite{K}, \cite{KR1}.
  Among them the unit sphere is the only compact Einstein space with
  $f'(t) = \sin t$ with the two zeros added (it is not a global warped product).
  
\item The case $a=0$ and $b=1$ leads to $f' = \gamma h'$ with a constant
  $\gamma$, the classical solutions of $(\nabla^2f)^\circ = 0$,
  a conformal gradient field \cite{K}.
\item Conversally, the assumption $f' = \gamma h'$ with a constant $\gamma$
  leads to $(\nabla^2f)^\circ = 0$ and thus to an Einstein space if $a \not = 0$.

    \end{enumerate}
  
  \end{pro}

\medskip

\begin{rem}
  A solution with $a = b = 1$ and a constant $c$ is normally called
  a {\it gradient Ricci soliton}, and the tensor
${\rm Ric}^\infty_f = {\rm Ric} + \nabla^2 f$ is called the 
{\it Bakry-\'Emery-Ricci tensor} \cite{WW}, for ${\rm Ric}^m_f
={\rm Ric}+\nabla^2 f-m^{-1} df \otimes df$ see
\cite{GD}.
A solution with $a = b = 1$ and a function $c$ is called a
{\it gradient Ricci almost soliton} \cite{PRRS}.
\end{rem}

\medskip
\begin{pro}
\
  {\rm (The case of low dimensions)}

  \begin{enumerate}
\item Any $2$-dimensional gradient Ricci soliton with $b \not = 0$
  is a warped product outside the set of critical points which are isolated.
  This follows from the known solutions of the equation $\nabla^2 f = \phi g$
  with any scalar function $\phi$ \cite{Ta}, \cite{K}, compare \cite{PW}.
\item Conversally, a $2$-dimensional warped product metric
  as above is a gradient
Ricci soliton if it admits a function $f(t)$ such that
$\nabla^2 f = -(K + c)g$ where $K$ denotes the Gaussian curvature
$K = -h'''/h'$. This implies $f'' = f'h''/h' = h'''/h' - c$ and,
consequently,  $f' = \gamma h'$ with a constant $\gamma$.
Locally this ODE has solutions for any $c$.

\item For $n=3$ the warped product approach with the equations above
  requires the extra assumption
that $M_*$ is a space of constant Gaussian curvature $\lambda_*$.
\item For $n=4$ the warped product approach with the equations above
  implies that $M_*$ is a space of constant
sectional curvature $\lambda_*$.
\end{enumerate}
\end{pro}

  \begin{exas} \label{cigar}
  1.  The metric
  $ds^2=dt^2+\tanh^2t \,du^2$ on $(t,u)\in [0,\infty)\times S^1$ 
  is a well known example of a $2$-dimensional
    complete Ricci soliton defining a 
    complete metric on $\R^2,$ 
     cf. K\"uhnel~\cite[Example 4(a)]{K}. Here     
    $c=0, \gamma=-2, h'(t) = \tanh t, f' = -2h', h''(t) = (\cosh t)^{-2}, K = -h'''/h' = 2 h'' = -f'', \nabla^2 f = f'' g = -2h'' g=K \,g.$ 
    With $x=\sinh t \, \sin u; y=\sinh t \cos u$
    this metric on $\R^2$ is
    of the form
    $$ds^2=\frac{dx^2+dy^2}{1+x^2+y^2}; f(x,y)=-\log (1+x^2+y^2),$$
   it is called \emph{Hamilton's cigar soliton} 
   in the literature, cf.~\cite{Hamilton1988}.

    2. The product of an Euclidean $\mathbb{R}^m$ with a $k$-dimensional
    Einstein space with Einstein constant $\lambda \not = 0$
    is a gradient Ricci soliton for the function
    $f(x) = \lambda |x|^2$ and for $c = -\lambda$.

  \end{exas}

\section{Complete and explicit solutions} 
\label{sec:complete}

Throughout this section we consider Ricci almost solitons 
in a slightly generalized form as solutions of
$(a\rm{Ric} + b \nabla^2 f)^\circ = 0$ with given constants $a$ and $b$.

Complete solutions are of particular interest.
Proposition~\ref{given-h} gives a possibility
for constructing Ricci almost solitons as warped products.
In fact there are complete Ricci almost solitons
that can be explicitly described in terms of elementary functions.
For $n=2$ the situation is better because one can prescribe the three
constants $a, b, c$ and solve the resulting ODE for one function
(not coupled as in the case $n \geq 3$).

\begin{pro}
  If $n=2$ 
  then the metric $g_*$ is $1$-dimensional,  hence $\lambda_*=0,$
  and the two coupled equations
  \eqref{eq:t-direction} und \eqref{Hauptgleichung} read as:
 \begin{eqnarray*}
- a h''' + b h'f'' + c h' &=& 0\\
a h'h''' - bf'h'h'' - ch'^2 &=& 0
 \end{eqnarray*}
For $h' \not = 0$ these are equivalent to the following two uncoupled equations:
 \begin{eqnarray}
f' &=& \gamma h' \\
\label{eq:ah3}
- a h''' + b \gamma h'h'' + c h' &=& 0
 \end{eqnarray}
for some constant $\gamma \not = 0$ (the case $\gamma = 0$ is trivial).
The second equation allows explicit solutions at least for the
inverse function $t = t(h')$.
\end{pro}

{\it Proof.}
The sum of both equations (the second divided by $h'$) leads to 
$b(h'f'' - h''f') = 0$, so we have either  b=0  or
$\log(f')' = \log(h')'$ which is equivalent to $f' = \gamma h'$ for a constant $\gamma.$ 
Then both equations are converted into (11). 
For $a=1, b=0$ see Prop. 2.2 (1).
\hfill $\Box$

\medskip
There is a method for solving an ODE of type $y'' = F(y,y')$
\cite[\S 11]{Wa}. We concentrate on the case $a = b = 1$.

\begin{lem}
  For $c=0$ Equation~\eqref{eq:ah3} is of first order by
  $(h'' - \frac{\gamma}{2}h'^2)' = 0$. There are various odd functions
  $h'$ solving this equation.
  Typical solutions with initial condition $h'(0) = 0$ are
  $h'(t) = \tan t$ for $\gamma = 2$ and
  $h'(t) = \tanh t = \frac{e^{2t}-1}{e^{2t}+1}$ for
  $\gamma = -2$ (the complete \emph{cigar metric}.)
  Other typical solutions are $h'(t) = \frac{1}{t}$ for $\gamma = -2$ and
  $h'(t) = \coth t = \frac{e^{2t}+1}{e^{2t}-1}$ for $\gamma =-2$.
  There are no even functions
  solving this equation.
\end{lem}

{\it Proof.} In the first case we have $h'' - h'^2 = 1$, in the
second case $h'' + h'^2 = 1$, in the third case $h'' + h'^2 = 0$,
in the fourth case $h'' + h'^2 = 1$.

\begin{exa}
  Let $h(t) = \arctan(t)$, then for any $\lambda_*$
  there is a global function $f(t)$ such that $\mathbb{R} \times_{h'} M_*$
  is a complete Ricci almost soliton whenever $M_*$ is complete
  with Einstein constant $\lambda_*$.
  \end{exa}

{\it Calculation.} We have $h'(t) =\frac{1}{t^2 + 1} > 0$ and
$h''(t) = -\frac{2t}{(t^2 + 1)^2}$. For any constants
$a,b \not = 0$ we have to solve Equation~\eqref{neueGleichung}:
$$a(n-2)\big[h''^2 - h'h''' \big]
 + bh'^3 \Big( \frac{f'}{h'}\Big)' 
=-a(n-2)\frac{2(t^2 - 1)}{(t^2 + 1)^4} + b \frac{1}{(t^2 + 1)^3} \big((t^2 + 1)f'\big)' = a \lambda_*.$$
This leads to
$$\big((t^2 + 1)f'\big)' = \frac{a \lambda_*}{b}(t^2 + 1)^3 + \frac{a (n-2)}{b}\cdot\frac{2(t^2 - 1)}{t^2 + 1}$$ and

$$(t^2 + 1)f' = \frac{a \lambda_*}{b}\int(t^2 + 1)^3 dt + \frac{a (n-2)}{b} \int \frac{2(t^2 - 1)}{t^2 + 1} dt$$
$$= \frac{a \lambda_*}{b}\Big(\frac{t^7}{7} + \frac{3t^5}{5} + t^3 + t \Big)
+ \frac{ a (n-2)}{b}\Big( 2t - 4\arctan(t) + C\Big).$$
For obtaining $f(t)$ itself it only remains to evaluate the integral in
\begin{eqnarray*}
f(t)& =& \frac{a \lambda_*}{b} \int \Big(\frac{t^7}{7} + \frac{3t^5}{5} + t^3 + t \Big)\frac{dt}{t^2 + 1} \\
&&+ \frac{a (n-2)}{b} \Big(\log(t^2 + 1) - 2\big(\arctan(t)\big)^2 + C\arctan(t) + D \Big)
\end{eqnarray*}

with constants $C,D$.
All expressions are globally defined on
$\mathbb{R}$ and give a global solution with $h'(t) \not = 0$.

\medskip
If $g_*$ is the standard metric on the sphere, then the space can
be viewed as a rotationally symmetric deformed infinite cylinder with
a ``cusp'' at $\pm \infty$.
This satisfies the equation
$$a({\rm Ric})^\circ + b(\nabla^2 f)^\circ = 0\ \mbox{ or } \ 
a{\rm Ric} + b\nabla^2 f + cg = 0$$
with $c = -(aS + b\Delta f)/n$.

\bigskip

\begin{exa}
  Let  $h(t) = \frac{1}{3}t^3 + t$, then for any $\lambda_*$
  there is a global function $f(t)$ such that $\mathbb{R} \times_{h'} M_*$
  is a complete Ricci almost soliton whenever $M_*$ is complete
  with Einstein constant $\lambda_*$.
  \end{exa}

{\it Calculation.} We have $h'(t) = t^2 + 1 > 0$ and $ h''(t) = 2t$.
This implies
$h''^2 - h'h''' = 4t^2 - 2(t^2 + 1) = 2(t^2 - 1)$. Therefore we have
to solve the equation
$$\Big(\frac{f'}{t^2 + 1}\Big)' = \Big(\lambda_*a - 2a(n-2)(t^2 - 1)\Big)\frac{1}{b(t^2 + 1)^3},$$
and $f'(t)$ can be determined by integration of rational functions.
In a second step one can express $f(t)$ at least by integrals of rational
functions and their integrals.
Obviously $f(t)$ is defined for any $t \in \mathbb{R}$

\medskip
If $g_*$ is the standard metric on the sphere, then the space can
be viewed as a rotationally symmetric deformed hyperboloid which opens up
quadratically at the two ends.

\begin{exa}
  Let $h(t) = (t^2 + 2)e^t$, then for any $\lambda_*$
  there is a global function $f(t)$ such that $\mathbb{R} \times_{h'} M_*$
  is a complete Ricci almost soliton whenever $M_*$ is complete
  with Einstein constant $\lambda_*$.
  \end{exa}

{\it Calculation.}
We have $h'(t) = (t^2 + 2t + 2)e^t > 0, h''(t) = (t+2)^2 e^t,
h'''(t) = (t^2 + 6t + 8)e^t$. 
We have to solve the equation
$$a(n-2)\big[h''^2 - h'h''' \big]
 + bh'^3 \Big( \frac{f'}{h'}\Big)'$$ 
 $$=-2a(n-2)(t^2 + 12t + 8)e^{2t} + b (t^2 + 2t + 2)^3e^{3t} \Big(\frac{f'}{(t^2 + 2t + 2)e^t}\Big)' = a \lambda_*.$$
 $f'$ can be determined by explicit integration, $f$ can at least be expressed
 by integrals.

 \medskip
 If $g_*$ is the standard metric on the sphere, then the space can
be viewed as a rotationally symmetric deformed infinite \emph{trumpet} with
a \emph{cusp} at one end and a huge opening at the other.

\bigskip

\begin{rem}[The power series approach for Ricci solitons]

For $a=b=1$ and a constant $c$ one can try to solve the system of equations
\ref{eq:t-direction} und \ref{Hauptgleichung}
simultaneously by a power series approach:
$$h'(t) = \sum_{i\geq 0}a_it^i, \quad f'(t) = \sum_{j\geq 0}b_jt^j$$
For $c=-1, a_0 = b_0 = 0, a_1 = b_1 = 1$ and $\lambda_* = n-2$ we have the classical solution
(Gaussian Ricci soliton) with $h'(t) = f'(t) = t$.
This is nothing but the Euclidean space in polar coordinates with the
function $f(x) = \frac{1}{2}|x|^2$, compare Proposition~\ref{specialcases}(2).
Under the assumption that both $h'$ and $f'$ are odd functions,
one can transform Equation~\eqref{eq:t-direction} into 
a (very complicated) recursion formula for the coefficients
depending on given constants $a_1, b_1, c$ and
$\lambda_*$.
Note, however, that even for the complete \emph{cigar} solution in
Examples~\ref{cigar} (1)
the power series has a finite radius of convergence.
\end{rem}
 \section{Ricci-Bourguignon solitons}
 An $n$-dimensional {\it Ricci-Bourguignon soliton}
 is a self-similar solution of the flow with any tensor
 ${\rm Ric} - \rho Sg$ instead of the Ricci tensor, where $\rho$ is a
 constant \cite[3.24]{B}, \cite{D}.
 The case $\rho = \frac{1}{2(n-1)}$ is referred to also as a
{\it Ricci-Schouten soliton} \cite{Bo}.
 In particular we are going to deal with the 
traceless Ricci tensor instead of the Ricci tensor itself, so we have
$\rho = \frac{1}{n}$.
We are going to investigate solutions of the equation 
\begin{equation}\label{Bourguignon}
{\rm Ric} + \nabla^2f - \frac{S}{n} g = \lambda g
\end{equation}
with a constant $\lambda$ where $S$ denotes the scalar curvature \cite{D}.
The case $\lambda=0$ is then
referred to as a {\it steady gradient Ricci-Bourguignon soliton}.

\begin{rems} 1. The only complete Einstein spaces satisfying
    Equation~\eqref{Bourguignon}
    are the Euclidean space in polar coordinates
    (if $\nabla^2f  = \lambda g$ with a constant
    $\lambda \not \equiv 0$)
and products of an interval with a complete Ricci flat manifold
(if $\nabla^2 f = 0$).
This follows from the classification of Einstein spaces admitting a solution
of the equation $\nabla^2 f = \phi g$ with a scalar function $\phi$
\cite{K}, \cite{KR1}.

\smallskip
2. If the gradient of $f$ is assumed to be a conformal vector
field, then Equation~\eqref{Bourguignon} implies that $M$ is Einstein.
Consequently, even if $\lambda$ is admitted to be a non-constant function,
the only complete cases are listed in \cite[Thm.27]{K}, implicitly
already in \cite[Sect.3]{Ta}. In particular,
the only compact case is the standard sphere. Implicitly this result can 
already be found in \cite[Thm.1]{Ta} since an Einstein space conformal
to the standard sphere is isometric to the standard sphere. 
It is restated in \cite[Thm.1.5]{D}.
For the last statement a constant scalar curvature
is sufficient \cite[Thm.24]{K}.

\smallskip
3. Complete gradient Ricci-Bourguignon almost solitons can be constructed by the
same principle as in Section~\ref{sec:complete},
compare \cite[Ex.1.7]{D}.

\smallskip
4. The case $n=2$ is not interesting here since the traceless Ricci tensor
always vanishes.
\end{rems}
\begin{pro} \label{steady}
  A warped product $M = \mathbb{R} \times_{h'} M_*$ of dimension $n \geq 4$
  with a potential function
  $f(t)$ is a steady gradient Ricci-Bourguignon soliton if and only
  if $f$ is harmonic and Equation~\eqref{alteGleichung} is satisfied with
  $a=b=1$.

  The same holds for $n=3$ under the assumption that $g_*$ is a metric
  of constant curvature $\lambda_*$.
  \end{pro}
    {\it Proof.}
    From the trace of Equation~\eqref{Bourguignon} we obtain
    $\Delta f = \lambda = 0$.
    Furthermore Equation~\eqref{alteGleichung} holds in general for
    this type of warped products, as shown in Section 2.

    Conversely, Equation~\eqref{alteGleichung}
    with $a=b=1$ implies ${\rm Ric} + \nabla^2f + c g = 0$
    for some function $c$. From the trace of this equation and $\Delta f = 0$
    one sees that $c = -\frac{S}{n}$. \hfill $\Box$
    
\begin{cor}\label{twoequations}
\qquad   An $n$-dimensional warped product metric \newline
  $ds^2 = dt^2 + (h'(t))^2ds_*^2$
  is a steady gradient Ricci-Bourguignon soliton if and only if
  $h$ and $f$ satisfy the following equations with the constant $\lambda_*$
  as the normalized scalar curvature of $g_*$:
\begin{eqnarray}
 (n-2)\big[h''^2 - h'h'''\big] 
 + h'\big[f''h' - f'h'' \big] 
&=& \lambda_* \\
\Delta f = f'' + (n-1) \frac{f'h''}{h'} &=& 0.
\end{eqnarray}
\end{cor}

\begin{thm} \label{thm:Ricci-Bourguignon}
  All complete warped products in dimension $n \geq 3$
  according to Proposition~\ref{steady}
and with $\lambda_* = 0$ can be classified as follows:

$(M_*,g_*)$ must be complete and Ricci flat, the possible
non-constant functions $h' > 0$ and $f'$ are the following:
$$\textstyle h'(t) = \Big(\big[\frac{B}{D} + E^{n-1}\big]e^{(n-1)Dt} - \frac{B}{D}\Big)^{1/(n-1)} \ \ \mbox{ and } \ \
f' = C (h')^{1-n}$$
with constants $B > 0, C > 0, D<0, E = h'(0)$ such that
$E^{n-1} > - \frac{B}{D}$.
\end{thm}

{\it Proof.} The case of an Einstein space can be ruled out here since
there are only the examples with a constant $h'$, compare the Remark 1 above.
The Euclidean space in polar coordinates with $h' = t$ and $\lambda_* = 1$
is not an example because $\nabla^2f = 0$
is not possible in this case for a function $f(t)$
although there is a harmonic function $f(t)$.
Now we assume that $h'$ is not constant and is nowhere zero.
From Corollary~\ref{twoequations} we obtain that the two equations
\begin{eqnarray*}
 (n-2)\big[h''^2 - h'h'''\big] 
 + h'\big[f''h' - f'h'' \big] 
&=& 0 \\
f'' + (n-1) \frac{f'h''}{h'} &=& 0.
\end{eqnarray*}
are satisfied.
The second equation is equivalent to
$$\textstyle\frac{f''}{f'} = -(n-1)\frac{h''}{h'} \ \mbox{ or } \ \log(f') = -(n-1)\log(h') + C_0$$ with a
  constant $C_0$, so it is equivalent to
\begin{equation}
  f' = C (h')^{1-n}
  \end{equation}
    with a constant $C>0$.
  Therefore $f$ is globally defined whenever $h'$ is globally defined without
  a zero. In a zero of $h'$ the function $f'$ tends to infinity, so
  it is really a singularity.
  This case will occur, depending on the choice of an additional
  constant $D$ that comes in. 

\medskip
  If we insert the second into the first equation we obtain
  \begin{equation*}
 (n-2)\big[h''^2 - h'h'''\big] 
 + h'\big[- n f'h'' \big] 
= 0
  \end{equation*}
  or, equivalently
\begin{equation}\label{Gleichung}
\textstyle    (2-n) \Big( \frac{h''}{h'}\Big)'
 + \big[- nC h''(h')^{-n} \big] 
= 0.
  \end{equation}

  \medskip
  We can solve Equation~\eqref{Gleichung} by
$$\textstyle h''(h')^{-n} = \frac{1}{1-n} \Big((h')^{1-n}  \Big)'$$
which implies
$$\textstyle(2-n) \frac{h''}{h'} = \frac{nC}{1-n}(h')^{1-n} + D(2-n)$$
with another constant $D$.
Now the discussion splits into the two cases
$$h'' = B (h')^{2-n}$$
with the constant $B = \frac{nC}{(n-1)(n-2)}>0$  \ if $D=0$ \ and
$$h'' = B (h')^{2-n} + Dh'$$ otherwise.

\bigskip
\underline{1. The case $D=0$}

Here the solution is  $$h'(t) = \Big( B (n-1) t + E^{n-1} \Big)^\frac{1}{n-1}$$
with an additional constant $E$ according to the initial condition
$h'(0) = E$. Necessarily $h'$ has a zero at $t=t_0$ for a certain $t_0$.
Therefore this solution does not lead to a complete
manifold. The singularity at $t=t_0$ is recognized by
an infinite value of $h''$ and also of $f'= C\big(h'\big)^{1-n}$.
The qualitative behavior does not depend on the choice of the constant $E$.

\bigskip
\underline{2. The case $D\not = 0$}

The general solution is
$$\textstyle h'(t) = \Big(\big[\frac{B}{D} + E^{n-1}\big]e^{(n-1)Dt} - \frac{B}{D}\Big)^{1/(n-1)}$$
with an additional constant $E$ according to the initial condition $h'(0) = E$.
Obviously this solution $h'(t)$ is globally defined and strictly positive
if and only if $D<0$ and $E^{n-1} > -\frac{B}{D} > 0$. Recall $B>0$ in any case.
In the limit $t \to \infty$ the metric is asymptotic to the
product metric $dt^2 + \big|\frac{B}{D}\big|^{2/(n-1)}ds_*^2$
on $\mathbb{R} \times M_*$ which is Ricci flat.
Strictly negative solutions can be considered only if $n$ is even but these
cases don't give any really new solutions for the warped product metric.
\hfill $\Box$

\medskip
    {\bf Remarks.} 1. For fixed $B$ and $E$ the transition
    $D \to 0$ in the second case
    leads to the limits 
    $\lim_{D \to 0} B\frac{e^{(n-1)Dt} - 1}{D} = B(n-1)t$ and
    $\lim_{D \to 0} e^{(n-1)Dt} = 1$ for any $t$
    (pointwise), so the limit function coincides with the solution
    in the case $D=0$. The singularity (a zero of $h'$)
    occurs if $E^{n-1} < -\frac{B}{D}$, and $h'$ is constant
    for $E^{n-1} = -\frac{B}{D}$.

    \medskip
    2. The global shape in the case $D \not = 0$ can look like     
    a huge \emph{trumpet} with a compact Ricci flat space as the fibre,
    asymptotically to a cylinder at one end and an infinite opening at the other.
    In the case $D=0$ it looks like a \emph{cigar}, modelled by the function
    $h'(t) = \sqrt[n-1]{t}$. If $M_*$ is compact and Ricci-flat then there
    is a topological singularity at $t=0$.

\medskip
3. In the complete examples with $D \not = 0$ 
in Theorem~\ref{thm:Ricci-Bourguignon} 
one can compute the functions
$f, h'',h'''$ explicitely
with the positive constants $\alpha:=E^{n-1}+B/D>0, \beta:=-(n-1)D>0, \gamma:=-B/D>0,
C>0$ and control their signs. We obtain
the following equalities resp. inequalities for all $t \in \R:$
\begin{eqnarray*}
    f(t)&=&\frac{C}{\beta \gamma}\ln \left(\alpha e^{-\beta t}+\gamma\right)+\frac{C}{\beta} t\\
    h'(t)&=&\left(\alpha e^{-\beta t}+\gamma\right)^{\frac{1}{n-1}}>0;
    h''(t)<0;
    h'''(t)
    >0.
\end{eqnarray*}
Then the scalar curvature $R=R(t)$ is negative, since
by \cite[Lem.13]{K}
\begin{equation}
h'^2(t)R(t)=-\frac{n-2}{n} h''^2(t)-\frac{2}{n}h'(t)h'''(t)<0
\end{equation}
for all $t \in \R,$ and
hence $\lim_{t\to \infty}R(t)=0.$
Also note that by \cite[Lem.13]{K}:
$$\textrm{Ric}(\partial_t,\partial_t)=-(n-1)\frac{h'''(t)}{h'(t)}<0.$$
In  \cite[Thm. 4.3]{CM} for different values of 
$\rho$ the existence 
and uniqueness 
of Ricci-Bourguignon solitons (named
\emph{gradient steady $\rho$-Einstein solitons} in \cite{CM})
is shown on warped products
with $\textrm{Ric}(\partial_t,\partial_t)>0$
for some $t.$  

\medskip
4. The same method can be applied to the classification of rotationally symmetric Ricci-Bourguignon solitons only by solving the ODEs in Corollary 4.2 also for $\lambda_* > 0.$ We do not see a way to achieve that by direct and explicit integration. Presumably one would have to employ other methods.
The situation is similar in the case of Bryant's Ricci solitons in dimensions $n >2$: 
The solutions exist but not in terms of elementary functions.

\section*{Statements and Declarations}
The authors have no relevant financial or non-financial interests
to disclose.

\smallskip

Data sharing is not applicable to this article as no datasets were generated or analysed during the current study.  

\bigskip
{\scriptsize

\noindent 

Wolfgang K\"uhnel,
Fachbereich Mathematik, 
Universit\"at Stuttgart, 
70550 Stuttgart, Germany,\\
E-mail: {\tt kuehnel@mathematik.uni-stuttgart.de}

\smallskip

\noindent 

Hans-Bert Rademacher,
Mathematisches Institut,
Universit\"at Leipzig,
04081 Leipzig, Germany,\\
E-mail: {\tt rademacher@math.uni-leipzig.de}
}

\end{document}